\RequirePackage{fix-cm}
\documentclass[smallcondensed]{svjour3}     
\smartqed  
\usepackage{graphicx}
\usepackage[utf8]{inputenc}
\usepackage[english]{babel}
\usepackage{tikz,xcolor,hyperref}
\usepackage{marvosym}
\usepackage{amsmath}
\usepackage{enumerate}

\begin{document}
	
	\title{Dominance Move calculation using a MIP approach 
	}
	\subtitle{for comparison of multi and many-objective optimization solution sets}
	
	
	\author{
		{${\textbf{Claudio Lucio do Val Lopes}}^{1}$}   \and
		{${\textbf{Flávio Vinícius Cruzeiro Martins}}^{1}$}   \and
		{${\textbf{Elizabeth Fialho Wanner }}^{1,2}$} 
	}
	

	
	\institute{
		\Letter   $  $   Claudio Lucio do Val Lopes\\
		\href{claudiolucio@gmail.com}{claudiolucio@gmail.com} \\
		{$\textbf{ }^{1}$}CEFET-MG, Programa de pós-graduação em modelagem matemática e computacional. Av. Amazonas, 7675 - Belo Horizonte,MG, Brazil - 30510-000\\
		$\textbf{ }^{2}$ Computer Science Group, Aston University, Birmingham, UK\\   
	}
	
	\date{Received: date / Accepted: date}

	\maketitle
	
	\begin{abstract}
		Dominance move (DoM) is a binary quality indicator that can be used in multiobjective optimization. It can compare solution sets while representing some important features such as \textit{convergence}, \textit{spread}, \textit{uniformity}, and \textit{cardinality}. DoM has an intuitive concept and considers the minimum move of one set needed to weakly Pareto dominate the other set. Despite the aforementioned properties, DoM is hard to calculate. The original formulation presents an efficient and exact method to calculate it in a biobjective case only. This work presents a new approach to calculate and extend DoM to deal with three or more objectives. The idea is to use a mixed integer programming (MIP) approach to calculate DoM. Some initial experiments, in the biobjective space, were done to verify the model correctness. Furthermore, other experiments, using three, five, and ten objective functions were done to show how the model behaves in higher dimensional cases. Algorithms such as IBEA, MOEAD, NSGAIII, NSGAII, and SPEA2 were used to generate the solution sets, however any other algorithms could be used with DoM indicator. The results have confirmed the effectiveness of the MIP DoM in problems with more than three objective functions. Final notes, considerations, and future research are discussed to exploit some solution sets particularities and improve the model and its use for other situations.

	\end{abstract}

	\section{Introduction}
	\label{intro}
	Problems with conflicting objectives arise in most real world optimization cases. These problems are circumscribed by multi or many-objective optimization techniques \cite{Chugh:2019:SHC:3333833.3333883} and exist in different domains such as \cite{Gonzlezlvarez2013AnalysingTS}, \cite{Deb_Kalyam_2008}.  The solution sets for these problems are formed in such a way that each objective represents a trade-off among other competing objectives. 
	
	In contrast to the mono-objective optimization, the solutions generated by multiobjective algorithms can be difficult to obtain and compare. Such difficulty grows as the number of candidate solutions and the objectives increase. When there are two or three objectives only, some graphical techniques help to examine the solution set visually. However, when the number of objectives is greater than three, this task is challenging, needing more advanced visualization techniques that can present location, shape, and solution set distribution \cite{DBLP:journals/swevo/IbrahimRMD18}.
	
	Quality indicators are suitable for situations when we need to sum up the solution set \cite{1197687}, taking into account their characteristics. They have been used to compare the outcomes of multiobjective algorithms quantitatively (Section \ref{definitions}). In a recent paper \cite{Li:2019:QES:3320149.3300148}, 100 quality indicators were discussed considering some state of the art indicators focusing on which quality aspects these indicators have, as well as its strengths and weaknesses.
	
	In \cite{DBLP:journals/corr/LiY17a}, a new quality measure, called dominance move (DoM) was proposed. This measure is able to capture all aspects of solution sets' quality such as Pareto compliant, spread, uniformity, and cardinality (Section \ref{relatedwork}). DoM measures the minimum `effort' that one solution set has to make in trying to dominate another set, more specifically the sum of the movement (i.e., Minkowski distance) needed to make a set dominant. The authors presented an exact approach to calculate DoM for the bi-objective case. Despite the low computational cost, the DoM method cannot be applied directly or extended to problems with three or more objectives (Section \ref{DOM}).  
	
	DoM presents good quality aspects for an indicator, but its calculation in low computational cost is a challenge \cite{Li2015Thesis}, \cite{Li2015}, and \cite{Li:2019:QES:3320149.3300148}. This work focuses on a DoM formulation based on a mixed integer programming approach (MIP)\cite{DBLP:books/daglib/0023873} aiming to overcome this difficulty. The mixed integer formulation is presented, and initial experiments showed that the MIP DoM formulation is correct. More specifically, this paper presents the following contributions: \textbf{(i)} a MIP model to DoM calculation that must be valid for all dimensions and solution sets; \textbf{(ii)} the use of this model to calculate DoM for some common problem sets having three, five, and ten objective functions, including considerations that come from the optimal DoM solution (Section \ref{experiments}); \textbf{(iii)} questions and details about the model behaviour for some test sets, raising future research paths to tackle some problems (Section \ref{conclusion}).

	\section{Problems and definitions}\label{definitions}
	In general, a multiobjective optimization problem (MOP) includes \textbf{\textit{x}} decision variable vector from a decision space $\Omega \subseteq R^{N}$, and a set of \textit{M} objective functions. Without loss of generality, a minimization MOP can be simply defined as \cite{Yuan2018}:
	\begin{align}\label{MOP_definition}
	Min \quad F(x) = {[f_{1}(x), ... ,f_{|M|}(x)]}^{T},  &&  x \in \Omega 
	\end{align}
	
	The $F: \Omega \rightarrow \Theta \subseteq R^{M}$ is constructed by a set  of $M$ objective functions, which is a mapping from decision space $\Omega$ to vectors in \textit{M}-dimensional objective space $\Theta$. We are interested in the evaluation of these objective vector (solution) sets, and the comparison relation among them.
	
	Considering two solutions \textbf{p}, \textbf{q}  $\in \Theta$, it is possible to establish a relation in which \textbf{p} is said to \textit{weakly dominate} \textbf{q} if \textbf{$p_i$} $\leq$ \textbf{$q_i$} for 1 $\leq$ i $\leq$ M, it is denoted as  \textbf{p} $\preceq$ \textbf{q}. Moreover, if there is at least one objective $i$ on which \textbf{$p_i$} $<$ \textbf{$q_i$} then it is said that \textbf{p} \textit{dominates} \textbf{q}, and is denoted as \textbf{p} $\prec$ \textbf{q}.
	A solution \textbf{p} $\in \Theta$  is called Pareto optimal  if there
	is no \textbf{q}  $\in  \Theta$ that dominates \textbf{p}. The set of all Pareto optimal solutions of an MOP is called Pareto optimal frontier. In the same way, the weak dominance  relation can be stated to solution sets:
	\begin{definition}{\textit{Weak Dominance}\footnote{There are other common relations defined for solution sets, i.e. \textit{strictly dominance, dominance, better}, but for DoM context, the weak dominance concept is enough.}:} 
		The set \textbf{P} weakly dominates \textbf{Q}, i.e. it is denoted as \textbf{P}  $\preceq$ \textbf{Q}, if every solution q $\in $ Q is weakly dominated by at least one solution p $\in$ P .
	\end{definition}
	
	The comparison between two or more solution sets is of great importance. It can be used to compare the outcomes of multiobjective algorithms or even assist an algorithm during the search process for most proper candidates. It is paramount to make more precise statements in a quality indicator comparison, for example: if one algorithm is better than another, how much better is it? The following definition formalizes the quality indicators \cite{1197687}:
	
	\begin{definition}{\textit{Quality indicator}:} 
		An \textit{k}-ary quality indicator \textit{I} is a function \textit{I}:${\Theta}^k \rightarrow {\rm I\!R}$, which assigns each vector of \textit{k} solution sets {($P_1, P_2, ..., P_k$)} a real value I{($P_1$, $P_2$ , ..., $P_k$)}.
	\end{definition}

	The quality indicators can be unary, binary, or \textit{k}-ary, defining value to one solution set, two solution sets, or \textit{k} solution sets, respectively.
	In \cite{Li:2019:QES:3320149.3300148}, 100 indicators are listed and some are  discussed in detail. Four indicator quality facets are analyzed: \textit{convergence}, \textit{spread}, \textit{uniformity}, and \textit{cardinality}. These facets will be analysed and discussed in the following.

	
	It is expected that the Pareto dominance must be a central criterion in reflecting the \textit{convergence} of solution sets. For two solution sets $P$ and $Q$, for example, if $P$ weakly dominates $Q$ then $I(P, Q) = 0$. If P dominates some points of Q and Q does not dominate any point of P, it is reasonable to expect that the indicator prefer $P$ to $Q$.
	
	The \textit{spread} of one solution set must consider the region that the set is covering. It involves both the outer portion and inner portion of the set. It must be able to capture the extensity of solution sets accurately.
	
	The number of solutions in the set is another quality indicator facet named \textit{cardinality}.  Finally, a good indicator must prefer a set with uniformly distributed points, \textit{uniformity}, showing an equidistant spacing amongst solutions.
	
	It is plausible to add more details about the indicators due to its usage, one of them being the computational cost. Some indicators present all quality aspects but are hard to compute, notably in high dimensions or in high-cardinality sets. DoM and Hypervolume are well known examples. Other indicator details also deserve to be mentioned, such as the necessity for reference point or set, additional parameters, how to deal with scale, and normalization. 
	
	\section{Related work}\label{relatedwork}
	
	There are many indicators available, and they have been used in numerous situations in the literature \cite{Li:2019:QES:3320149.3300148}. Hypervolume (HV) \cite{8625504}, \cite{Yang2019},\cite{Bradford_2018}, and inverted generational distance(IGD) \cite{10.1007/978-3-319-15892-1_8} are some examples:
	\begin{itemize}{}{}
		\item \textbf{Hypervolume (HV)}:  Let $r^{'}= (r_1, ..., r_m)$ be reference points in the objective space that is dominated by all solution sets. Let $P$ be one solution set. The HV value of $P$ with regard to $r^{'}$ represents the volume of the region which is dominated by $P$ and dominates $r^{'}$. 
		
		\item \textbf{Inverted generational distance (IGD)}:  Let $R^{*}= ({r^{*}}_1, ..., {r^{*}}_m)$ be a reference set of uniformly distributed points on the Pareto front. Considering $P$ as an solution set to the Pareto front, the inverted generational distance between $R^{*}$ and $P$ is defined as:
		\begin{equation}\label{IGD}
		IGD(R^{*}, P) = \frac{ \sum\limits_{r \in R^{*}} {d(r,P)}}{|{R^{*}}|}
		\end{equation}
		
		Note that ${d(r,P)}$ is the minimum Euclidean distance from point $r$ to solution set $P$. This indicator is a measure that represents how far the solution set is from the Pareto front reference. Lower values of IGD represent a better performance. The IGD metric is able to measure both diversity and convergence of $P$ if $|{R^{*}}|$ is large enough \cite{cheng2018benchmark}.

	\end{itemize}
	There are other quality indicators such as Generational distance (GD) \cite{10.1007/978-3-319-15892-1_8}, $\epsilon$-indicator \cite{1197687}, KKTPM \cite{7283599}, and others \cite{Li:2019:QES:3320149.3300148} .
	
	An ideal quality indicator must present the four facets. Additionally, it must have a low computational cost, and it does not need a normalization (due to objective scale) and any additional parameters or reference points/sets.
	
	Dominance move (DoM) is an intuitive indicator, and it has the four desirable facets \cite{DBLP:journals/corr/LiY17a}. The first idea presenting DoM came from the performance comparison indicator, PCI \cite{Li2015}. Examining the  PCI proposal, it is quite similar to DoM in its essential purpose. PCI, a binary quality indicator, builds up a reference set using two solution sets, $P$ and $Q$. This reference set is then split up into clusters, and the indicator calculates the movement distance in order to weakly dominate the clusters. 
	
	Dominance move is a measure for comparing two sets of multi-dimensional points being classified as a binary indicator. It considers the movement of points in one set to make this set weakly dominated by the other set. The DoM does not need a priori problem knowledge, parameters, or a reference set. However, the computational cost to calculate DoM is prohibitive, and the method to obtain DoM for a number of objective functions higher than two is unknown. DoM can be defined as follows.
	
	\begin{definition}{DoM: }
		Consider that $P$ and $Q$ are sets of points, with  $ p_i$ points $ i \in \{{1,..,NP}\}$ and $q_j$ points $j \in \{{1,..,NQ}\} $. The dominance move of $P$ to $Q$, $DoM(P,Q)$, is the minimum total distance of moving points of P, such that any point in $Q$ is weakly dominated by at least one point in $P$. In fact, the problem is to find ${ \{{p^{'}}_1,{p^{'}}_2,..,{p^{'}}_{NP} \}}$ from  ${ \{p_1,p_2,..,p_{NP} \}}$ positions such that $P^{'}$ weakly dominates $Q$ and the total move from ${ \{p_1,p_2,..,p_{NP} \}}$ to ${ \{{p^{'}}_1,{p^{'}}_2,..,{p^{'}}_{NP} \}}$ must be minimum \cite{DBLP:journals/corr/LiY17a}.
	\end{definition}
	
	$P^{'}$ is a set of points that are candidates to dominate $Q$ with some update in one or many objectives that lead to a better distance such as expressed in Equation \eqref{DPQ}. In such way, it should be noted that $ p^{'}_k$ with $k \in \{{1,..,NP^{'}}\}$ must be generated from $p$.
	
	The formal expression of DoM can be stated as: 
	\begin{equation}\label{DPQ}
	DoM(P,Q) = \underset{P^{'} \preceq  Q}{minimize} \sum\limits_{i=1}^{NP} d(p_i,p^{'}_i) 
	\end{equation}
	
	The dominance move indicator is based on the properties of dominance relation among solutions that are trying to dominate each other. These solutions' efforts scale in a bottom up manner from the solutions to the set relations \cite{DBLP:journals/corr/LiY17a}. 
	
	The number of possibilities to find $P^{'}$ is numerous. Any combination of some $P^{'}$ can dominate $Q$. Consider $P^{'}_{(1..s)}$ from $P$ with some  movement updates in the objectives in conjunction with other $P^{'}_{(s+1..NP')}$ from original $P$, for example. 
	
	The number of such candidates is detailed in Equation (\ref{comb_plinha}), in which $g$ is a group with one or many $q_{j}$, and assuming that $p_{i}$ will be used as a base to be updated, generating $p^{'}_k$ that can weakly dominate all $g$ group. 
	\begin{equation}\label{comb_plinha}
	NP \sum\limits_{g=1}^{NQ}{NQ\choose g} =NP \left({NQ\choose 1} + {NQ\choose 2} + ... + {NQ\choose NQ}\right)= NP( 2^{NQ} - 1)
	\end{equation}
	
	The Equation \eqref{comb_plinha} is, in fact, an upper bound for the number of solutions. It is possible to discard some repeated $p^{'}_k$ generated from $g$ and $p_{i}$. It is worthy to note that the number of repeated candidates depends on the distributions for each solution sets in $g$  and $p_{i}$.
	
	The original DoM formulation proposes a method for the biobjective case \cite{Li:2019:QES:3320149.3300148}. The algorithm can be outlined as:
	\begin{list}{}{}
		\item \textit{Step 1:} Remove the dominated points in both $P$ and $Q$, separately. Remove the points of $Q$ that are dominated by at least one point in $P$.
		\item \textit{Step 2:} Denote  $R$ $=$ $P \cup Q$ and start the process. Each point of $Q$ in $R$ is considered as a group. For each group of $Q$, find its inward neighbor $r = n_{R}(q)$ in $R$. If the point $r \in P$, then merge $r$ into the group of $q$, otherwise $r \in Q$. If $r$ is not assigned to one group merge the two groups of $q$ and $r$ into one group.
		\item \textit{Step 3:} If there exists no point $q \in Q$ such that $q = n_{R}(n_{R}(q))$ (i.e., there is a loop between the points) in any group, then the procedure ends and there is an optimal solution to the case.
		\item \textit{Step 4:} There is a loop in one or some groups. The procedure replaces these loops with the ideal point. The ideal point is formed of the best of each objective in each point inside the loop or group. Return to step 3 until convergence.
	\end{list}
	The authors present all definitions, theorems and corollaries to prove that this algorithm is correct in the biobjective case. However, it is stated that there is no method for three or more objectives.
	
	\section{The MIP dominance move calculation approach} \label{DOM}
	
	Our goal is to modify the DoM formulation to deal with problems having three or more objective function. Our DoM calculation proposal is based on the perspective that the problem is, in fact, a particular instance of an assignment problem with two levels and some restrictions \cite{RePEc:eee:ejores:v:176:y:2007:i:2:p:774-793}. It is considered that to treat the problem, we have to find an assignment of $P$ to $Q$ with the restrictions that each $q_{i}$ must be assigned to one $p_{i}$ with the minimum distance. Nevertheless, in classic assignment problems, it is not possible that points  $p's$ change its own features, such as changing positions to alter the distances, for example. In DoM calculation, this issue must be considered. 
	
	A simple example to clarify the situation is given in \cite{Li:2019:QES:3320149.3300148}: consider $P$ as  \{(2.0, 2.0, 2.0), (2.0, 2.2, 1.5), (3.0, 1.6, 1.6)\}  and Q as \{(2.0, 1.2, 2.1), (2.0, 2.1, 1.0), (4.0, 1.5, 1.5)\}. The inward neighbor $r = n_{R}(q)$ of points $q_{1},q_{2}, $ and $q_{3}$ are, respectively, $p_{1},p_{2}, $ and $p_{3}$. This create an assignment of $P$ to $Q$ with the minimum $DoM(P,Q)$. Considering that $P$ is fixed: $DoM(P,Q) = d(p_{1},q_{1}) + d(p_{2},q_{2}) + d(p_{3},q_{3}) = 1.6$. However, if we considered a movement from $P$ to $P^{'}$, then $p_2$ would be transformed into $p^{'}_2 = (2.0,1.5,1.5)$. In this sense, we can find a better assignment and a lower value of $DoM(P,Q) = 1.5$. Clearly, other assignments from $P$ to $P^{'}$ and to ${Q}$ are capable to generate the same value (in \cite{Li:2019:QES:3320149.3300148}, the authors presented another answer with $DoM(P,Q) = 1.5$.)  
	
	\newcommand{\vc}[3]{\overset{#2}{\underset{#3}{#1}}}
	
	The problem can be modeled as a mixed integer programming approach. Generally, a MIP model (presented in Equation \ref{mip_general}) can be described as a set of variables $x_c$, $ c \in C$ that are non negative, and variables $x_i$, $ i \in I$ that are integer \cite{articleKlotz}. Additionally, $\vc{c^{T}}{}{C}$ and $A_{C}$ are the objective function and left-hand side constraint coefficients, respectively. It is important to note that they are non negative. In the same way, there are the $\vc{c^{T}}{}{I}$ and $A_{I}$ for integer variables. Lastly, for the constraint set, there is a $b$ right-hand side constant \cite{articleKlotz}.
	\begin{alignat}{2}\label{mip_general}
	&  \underset{}{\text{minimize }}
	\vc{c^{T}}{}{C}x_{C} + \vc{c^{T}}{}{I}x_{I}  \nonumber\\
	& \text{subject to}\nonumber \nonumber \\
	& A_{C}x_{C} + A_{I}x_{I} = b  \nonumber\\
	& x_{C}, x_{I} \ge 0   \nonumber\\
	&  x_{I} \in Integer   \\
	\nonumber
	\end{alignat}
	
	Using the MIP approach and DoM definition, $P$ and $Q$ are sets of points, with  $ p_i$ points $ i \in \{{1,..,NP}\}$ and $q_j$ points $j \in \{{1,..,NQ}\} $. $P$ and $Q$ are given in the problem. Otherwise, $P^{'}$ is a set of points, with  $p^{'}_k$ with $k \in \{{1,..,NP^{'}}\}$, which are candidates generated from $p_i$ to weakly dominate some $q_j$ (i.e., or a $g$ group), with some update in one or many objectives. It is relevant to note that if $p_i$ already dominates $q_j$ then $p^{'}_k = p_i$. $P^{'}$ weakly dominates  $Q$, resulting in a better distance such as expressed in Equation \eqref{DPQ}.
	
	The proposal aims to calculate the distances $d(p_i,p^{'}_k) $ and $d(p^{'}_k,q_{j})$ in the model and provide some limits for each $p^{'}_k$ that could vary due to $g$. Each $g$ group is formed by one or many $q_{j}$'s to be dominated. 
	\begin{alignat}{2}\label{mip_dom_model}
	&  \underset{}{\text{minimize}}
	\sum\limits_{i=1}^{NP} \sum\limits_{m=1}^{M} zp_{(i,m)}   +  \sum\limits_{i=1}^{NP} \sum\limits_{j=1}^{NQ} \sum\limits_{m=1}^{M} zpq_{(i,j,m)}  \\
	& \text{subject to}\nonumber  \\
	& zp_{(i,m)} \geq  0   \nonumber\\
	& zp_{(i,m)} \geq  p_{(i,m)} xp_{(i)}  - {p^{'}}_{(i,m)}   \nonumber\\
	& zp_{(i,m)} \leq  p_{(i,m)} xp_{(i)} &&  \forall \mathit{i} \in {NP}, \forall \mathit{m} \in {M} \label{restriction_zpio}\\
	& \nonumber\\
	& zpq_{(i,j,m)} \geq  0   \nonumber\\
	& zpq_{(i,j,m)} \geq  {p^{'}}_{(i,m)} - q_{(j,m)} -  p_{(i,m)} (1 - xpq_{(i,j)})    \nonumber\\
	& M_{(i,j,m)} = Max(0, p_{(i,m)} - q_{(j,m)})    \nonumber\\
	& zpq_{(i,j,m)} \leq  {p^{'}}_{(i,m)} - q_{(j,m)} + && \nonumber\\ 
	&(M_{(i,j,m)} - lbp_{(i,m)} - q_{(j,m)}) (1-xpqd_{(i,j,m)})\nonumber\\
	&zpq_{(i,j,m)} \leq  M_{(i,j,m)} (xpqd_{(i,j,m)})  && \forall \mathit{i} \in {NP},\forall \mathit{j} \in {NQ}, \forall \mathit{m} \in {M}
	\label{restriction_zpqijo}\\
	& \nonumber\\
	& ubp_{(i,m)} = p_{(i,m)}
	\nonumber\\
	& lbp_{(i,m)} = Min(p_{(i,m)}, Min(q_{(1..NQ, m)})) 
	\nonumber\\	
	& lbp_{(i,m)} \leq  {p^{'}}_{(i,m)} \leq  ubp_{(i,m)} && \forall \mathit{i} \in {NP},\forall \mathit{j} \in {NQ}, \forall \mathit{m} \in {M}
	\label{restriction_lowerandupper}\\
	%
	& \nonumber\\
	& xp_{(i)} \geq  xpq_{(i,j)}  ,   \label{restriction_xp}\\
	& xp_{(i)} \leq  \sum\limits_{j=1}^{NQ} xpq_{(i,j)}\label{restriction_sumxp} \\
	& \sum\limits_{i=1}^{NP}xpq_{(i,j)} = 1 ,  && \nonumber  \\
	%
	%
	& xp_{(i)} \in \{0,1\}, &&  \nonumber \\
	& xpq_{(i,j)}\in \{0,1\},  && \nonumber\\
	& xpqd_{(i,j,m)} \in \{0,1\}, &&  \forall \mathit{i} \in {NP}, \forall \mathit{j} \in {NQ} , \forall \mathit{m} \in {M}\label{restriction_xpodB}
	\end{alignat}		
	
	In the model expressed in Equation (\ref{mip_dom_model}) to Equation (\ref{restriction_xpodB}), these calculations can be done using MIP model. For each solution vector $p_i \in$ in $P$ and $q_j \in$ $Q$,  there are  $m \in \{{1,..,M}\}$ objectives. DoM was calculated using the auxiliary variables,  $zp_{(i,m)}$ that are equal to $d(p_i,p^{'}_k)$, and $zpq_{(i,j,m)}$ that are equal to $d(p^{'}_k,q_{j})$.
	
	The variable expressed in Equation (\ref{mip_dom_model}), $zp_{(i,m)}$ represents the distances between $p_{(i,m)}$ and a possible ${p^{'}}_{(i,m)}$ candidate. The constraints related to $zp_{(i,m)}$ are shown in Equation (\ref{restriction_zpio}). Essentially, it should be greater than zero and less than $p_{(i,m)}$ (always positive) and, finally, it must be greater than the difference in $p_{(i,m)}$ and ${p^{'}}_{(i,m)}$. 
	
	The variable $zpq_{(i,j,m)}$ represents the difference for each ${p^{'}}_{(i,m)}$ in an attempt to dominate $q_{(j,m)}$. The variable constraints are shown in Equation (\ref{restriction_zpqijo}). It ensures that ${p^{'}}_{(i,m)}$ can dominate $q_{(j,m)}$ using a binary variable $xpq_{(i,j)}$. Simultaneously, it will receive $Max(0,{p^{'}}_{(i,m)} - q_{(j,m)})$, as the difference to dominate $q_{(j,m)}$ or $0$. We have used a linearization technique to obtain the maximum function. Thus  $xpqd_{(i,j,m)}$, which assumes 1 to guarantee the maximum value, was introduced; otherwise, the solution will be unfeasible. It is worth mentioning that there is a region in which ${p^{'}}_{(i,m)}$ can assume values, an upper and lower bound, previously defined. 
	
	The $xp_{(i)}$, a binary variable in Equation (\ref{restriction_xpodB}), was used to guarantee that if a ${p^{'}}_{(i,m)}$ is used  in the model, we will know exactly which ${p_{(i,m)}}$ has generated the candidate. In the same way, the  $xpq_{(i,j)}$ reflects that a  ${p_{(i,m)}}$ generated a ${p^{'}}_{(i,m)}$  and is trying to dominate a $q_{(j,m)}$. These two binary variables were used to guarantee that at least one $q_{(j)}$ is going to be dominated by ${p^{'}}_{(i)}$. The constraint represented in Equation (\ref{restriction_xp}) is derived from the conjunctive normal form as described in Equation (\ref{restric_con}).
	
	\begin{alignat}{2}\label{restric_con}
	& \bigvee\limits_{j=1}^{NQ} xpq_{(i,j)}  \Rightarrow xp_{(i)} &&\nonumber \\  
	&\sim \bigvee\limits_{j=1}^{NQ}{\left(xpq_{(i,j)}\right)} \vee xp_{(i)} &&\nonumber   \\
	&\bigvee\limits_{j=1}^{NQ}{\left( \sim xpq_{(i,j)} \right)} \vee xp_{(i)} &&\nonumber   \\
	&\bigvee\limits_{j=1}^{NQ}{\left( \sim xpq_{(i,j)}  \vee xp_{(i)} \right)}  &&\nonumber   \\
	&\left(1- xpq_{(i,j)} +  xp_{(i)} \right) \geq 1  &&\nonumber\\
	&xp_{(i)}  \geq xpq_{(i,j)} , && \quad \quad \quad \quad \quad \quad \quad \forall_{(i,j)} \in I_{NP} \times I_{NQ}   \\
	\nonumber
	\end{alignat}

	In the same way, the constraints in Equation (\ref{restriction_sumxp}) were obtained from Equation (\ref{restric_sum_con}). All the $x's$ variables are binary and the constraints are expressed in Equation (\ref{restriction_xpodB}).
	
	\begin{alignat}{2}\label{restric_sum_con}
	& xp_{(i)} \Rightarrow  \bigvee\limits_{j=1}^{NQ} xpq_{(i,j)} &&\nonumber \\  
	&\sim xp_{(i)} \vee \bigvee\limits_{j=1}^{NQ}{\left(xpq_{(i,j)}\right)} &&\nonumber   \\
	&\left(1- xp_{(i)}\right) + \sum\limits_{j=1}^{NQ}{xpq_{(i,j)}}    \geq 1  &&\nonumber\\
	&xp_{(i)} \leq \sum\limits_{j=1}^{NQ}{xpq_{(i,j)}}    &&\quad \quad \quad \quad \quad \quad \quad \quad \quad \quad \quad \quad \quad \quad \\
	\nonumber
	\end{alignat}

	\section{Experiments}\label{experiments}
	
	\subsection{Biobjective case}
	
	The first test was done to show how MIP DoM calculation addresses the quality indicator facets: convergence, spread, uniformity, and cardinality  \cite{Li:2019:QES:3320149.3300148}. The same artificial experiments proposed in \cite{DBLP:journals/corr/LiY17a} to solve DoM in the biobjective case were adopted to verify the model correctness \footnote{The data were provided by the author, Dr Miqing Li.}.
	
	\begin{figure}[!htb]
		\centering
		{\includegraphics[width=\linewidth]{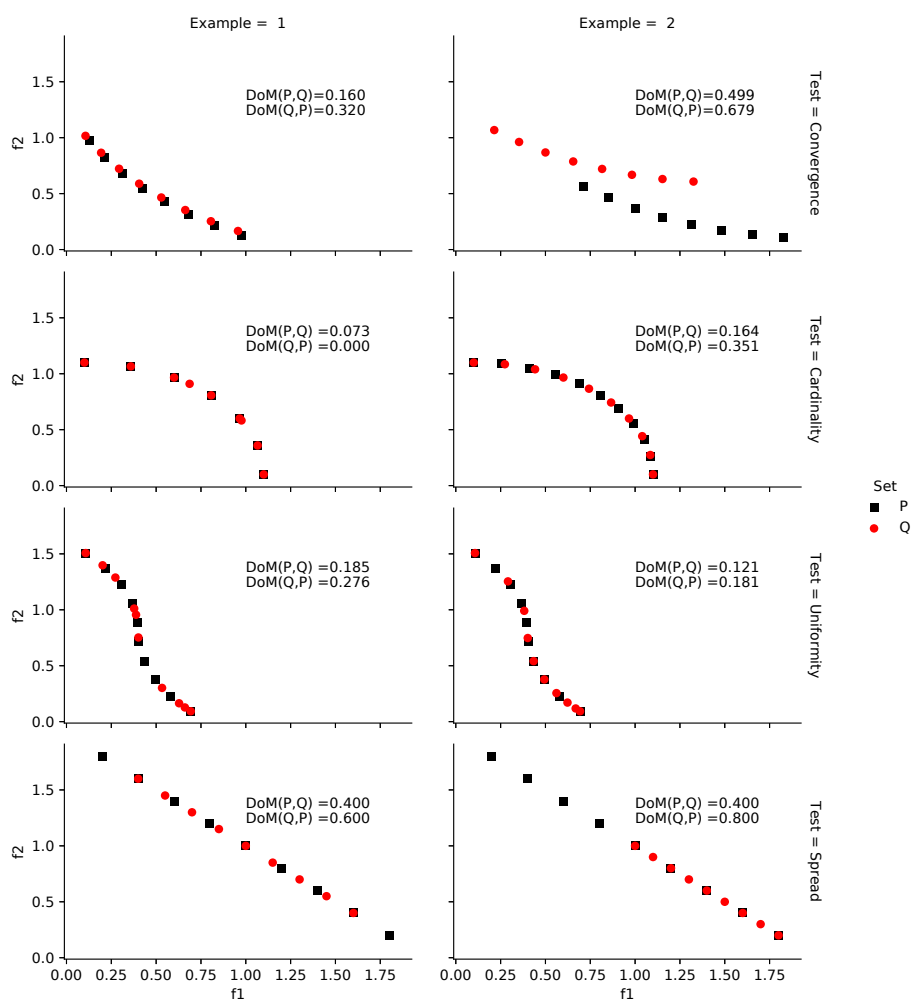} }
		\caption[ArtificialExperiments]{ Artificial experiments proposed in \cite{DBLP:journals/corr/LiY17a}  to assess the four facets of quality indicators: convergence, cardinality,  uniformity, and spread. For each facet, there was a combination of two graphics. In total, there were four rows of graphics: tests for convergence, test for cardinality,  test for uniformity, and test for spread. For each row, Examples 1 and 2 are given, and they are slightly different from one another. In each graphic there are two solution sets, \textit{P} and \textit{Q}.}
		\label{fig:ArtificialExperiments}%
	\end{figure}
	
	Figure \ref{fig:ArtificialExperiments} presents the artificial experiments for assessing the model correctness. Each facet of quality indicators, convergence, cardinality,  uniformity, and spread, was represented in a `row' in Figure \ref{fig:ArtificialExperiments}. Each row was composed of two graphics corresponding to two examples.  The examples are slightly different from one another. In each graphic there are two solution sets, \textit{P} and \textit{Q}.

	Convergence is an important factor to reflect Pareto dominance compliance of sets.  This behavior was shown considering two examples in Figure \ref{fig:ArtificialExperiments},  first row of graphics: test = convergence. In Examples 1 and 2, the $P$ and $Q$  are equals, however, in Example 1, the $P$ set was slightly improved in objective $f1$ , in comparison with $Q$. In Example 2, $Q$ has some points dominated by $P$, but not all of them. The MIP DoM values of both graphics reflected the dominance relation. 
	
	MIP DoM prefers solutions with more cardinality. In Figure \ref{fig:ArtificialExperiments}, test = cardinality, it could be observed that the solution sets had the same convergence, spread, and uniformity. In Example 1, $Q$ has one point more than $P$, and in Example 2, $P$ has two more points than $Q$.
	
	Uniformity indicates the preference for evenly distributed points. The solution sets in Figure \ref{fig:ArtificialExperiments}, test = uniformity, presented this feature. The sets have the same convergence, spread, and cardinality. Set $P$ was distributed uniformly, and $Q$ had a random distribution. In Example 1, $Q$ was distributed in the range of set P, and in Example 2, the distance between neighboring points in set $Q$ increased gradually from bottom to top. 
	
	Finally, MIP DoM must exhibit its preference for solutions with better spread. A set with better extensity has a smaller dominance move compared to its competitor. In Figure \ref{fig:ArtificialExperiments}, test = spread, considering Example 1, set Q was generated by shrinking P a little (more concentrated in the middle). In Example 2, set P was distributed uniformly in the range, while $Q$ assumed five bottom right points.
	
	The MIP DoM approach results were equal to the results found in \cite{DBLP:journals/corr/LiY17a}, using the biobjective algorithm proposed. The results agreement showed that the DoM MIP model is correct yielding to the same results obtained by its predecessor. 
	
	\subsection{Multi and Many-objective case}
	Following the initial and artificial test sets, it would be crucial to validate the MIP DoM model with more than two objectives. Firstly, a test set with three objectives was done and compared with classical indicators, such as HV and IGD, and also with visual graphics to assess the results. Secondly, an attempt to solve problems with five and ten objective functions was executed. In all tests, algorithms such as IBEA, MOEAD, NSGAIII, NSGAII, and SPEA2 were used to generate the solution sets. It is important to highlight that the goal was to assess the effectiveness of the proposed MIP DoM model and not to compare algorithm performance, so any other algorithms could have been applied to generate the solution sets. 
	
	
	For each problem set, one of the parameters chosen for each algorithm was the population size. The question is closely related to Equation (\ref{comb_plinha}) and the cardinality of solution set, one of the quality indicator facet. To have a good approximation set of the Pareto front, in terms of convergence, spread, and uniformity, the number of non-dominated solutions grows exponentially in relation to the dimension of the objective space. 
	
	In \cite{Sen1998MultipleCD}, the shape of the Pareto front was discussed in the niche sizes definition. A useful limit to the number of individuals in the population, given by $N = Mr^{M-1}$, was provided. The number of individuals in the population was $N$, $M$ was the number of objectives, and $r$ was the resolution required or the number of points needed to represent the Pareto front. This expression makes it clear how $NP$, or $NQ$, must be increased as the problem dimension grows, showing an exponential relation between $N$ and $M$. 
	
	However, other works in many-objective optimization did not strictly follow this rule. In \cite{8027123}, for example, the number of objectives was $M=\{3, 5, 10\}$ and the population size was $N = \{105, 126, 275\}$. In \cite{Yang2019}, an efficient hypervolume calculation was provided and some tests were done with $M=\{3, 4, 5\}$ and  $N = \{10, .., 200\}$. Likewise, in CEC'2018, a competition on many-objective optimization \cite{cheng2018benchmark} established   $M=\{5, 10, 15\}$ and the maximum population size as 240.
	
	Based on Equation (\ref{comb_plinha}), which was important in the proposed model (Equations (\ref{mip_dom_model})  to (\ref{restriction_xpodB})), we have decided to validate DoM using $M=\{3, 5, 10\}$ and $NP = NQ = \{50, 100, 200\}$ indicating the  final Pareto front size.
	
	All the experiments were done using \textit{Platypus} \cite{Brockhoff:2019:Platypus} and \textit{PyGMO} \cite{PyGMO} to generate the problem sets and to calculate the indicators (HV and IGD). The number of fitness evaluations was the same for all algorithms, 10000. The model in  Equation (\ref{mip_dom_model})  to (\ref{restriction_xpodB}) was implemented using \textit{Python} and \textit{GUROBI} \cite{gurobi} version 9.0.0 build v9.0.0rc2 running on a Linux 64 bits operational system with 12 CPU's and 16Gb of RAM.
	
	\subsubsection{Multiobjective cases}
	In \cite{Huband2011}, a review of multiobjective problem sets was presented. To validate our proposed approach, the aim was to select some well-known problem test sets in three dimensions, such as the \textit{DTLZ} and \textit{WFG} families. Some test sets were initially chosen based on the characteristics: convexity/concavity, disconnection, multimodality, and degeneracy. DTLZ1, DTLZ2, DTLZ3, and DTLZ7, which are linear, concave/multimodality, concave, and disconnected, respectively, were selected. In the same manner, WFG1, WFG2, WFG3 and WFG9 were also selected taking into account the properties: convex/mixed, convex/disconnected, linear/degenerate, and concave, respectively. 
	
\begin{figure}[!htb]
		{\includegraphics[width=\linewidth,height=10.5cm]{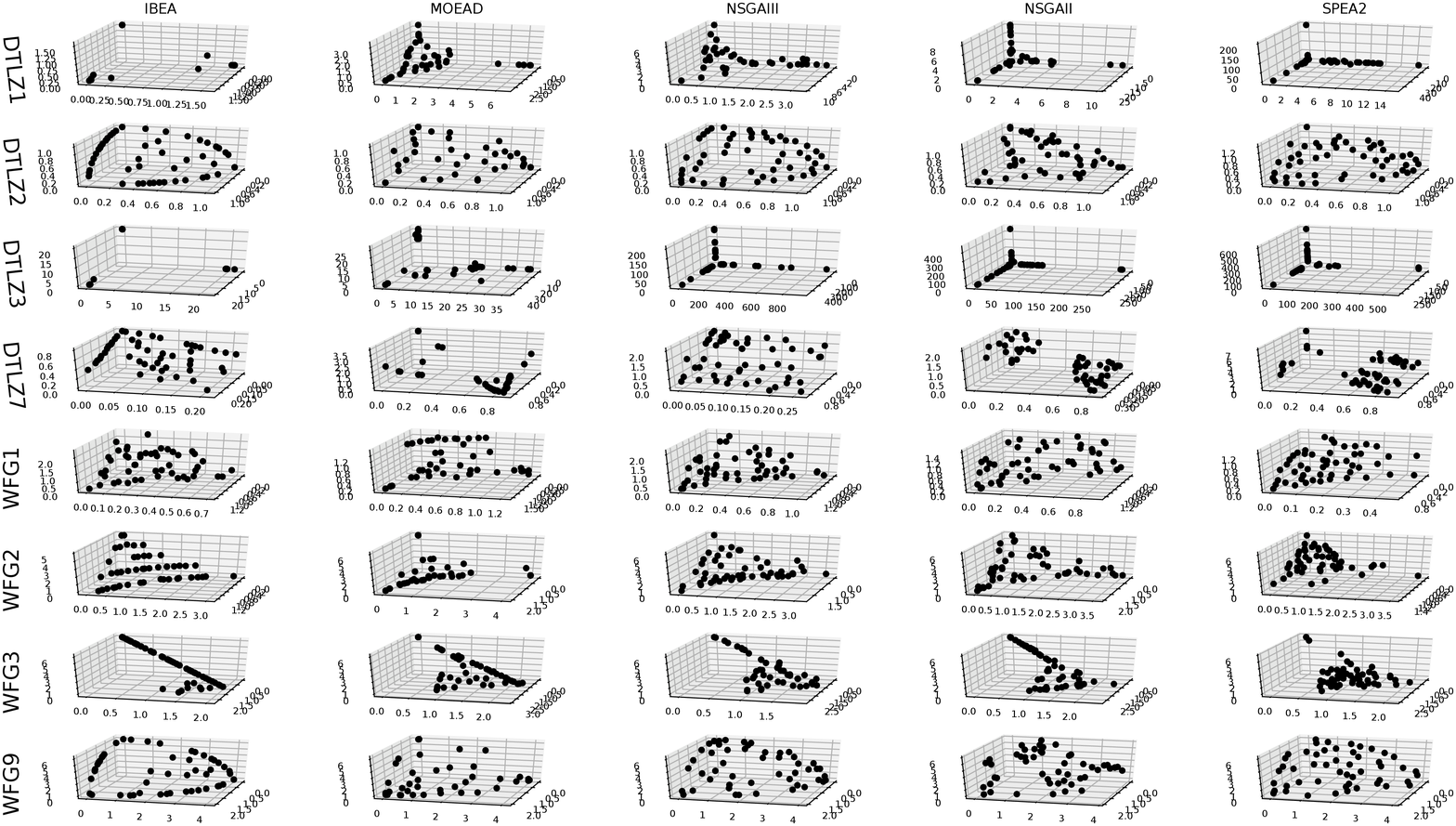} }
		\caption[WFG9]{Solution sets with \textit{NP} = \textit{NQ} = 50 solutions and three objectives, \textit{M} = 3, generated by IBEA, MOEAD, NSGAIII, NSGAII, and SPEA2 algorithms applied to DTLZ1, DTLZ2, DTLZ3, WFG1, WFG2, WFG3 and WFG9 problem sets.}
		\label{fig:WFG9_sample}%
	\end{figure}

	In parallel to the problem set, some algorithms were used to test and compare the MIP DoM with other indicators. Embedded with a similar purpose of that in\cite{ffb49e1a77a5478fabb03b092c89b2b7}, the NSGAIII and MOEAD algorithms (Pareto-based and decomposition approach, respectively)  were firstly chosen, and afterwards, IBEA \cite{Zitzler04indicator-basedselection}, SPEA2 \cite{Zitzler01spea2:improving}, and NSGAII \cite{Bradford_2018}. All the problem sets tested can be seen in Figure \ref{fig:WFG9_sample}.
	
	Tables \ref{tab:IGD} and \ref{tab:HV} show two unary quality indicators, the inverted generational distance (IGD) and hypervolume (HV). It is mandatory to have reference sets to calculate these indicators, and this task is not only a challenge one \cite{Li:2019:QES:3320149.3300148} but also sometimes is provided by the user \cite{Yang2019}. 
	We decided to use the maximum and minimum values, amongst all algorithm solutions and objectives, for HV and IGD.
	
	It is noted that DTLZ1 and DTLZ3, for example, had the HV value inflated by the presence of the dominance resistant solutions \cite{Li:2019:QES:3320149.3300148}, non-dominated solutions with a poor value in one objective but with good values in others. It is also important to observe that IGD and HV are highly sensitive to the reference set. If such reference set does not represent well the Pareto front, then the indicator is compromised. This parameter is crucial to these indicators' validity. 
	
	\begin{table}[!htb]
		\caption{IGD quality indicator for problem sets generated using DTLZ1, DTLZ2, DTLZ3, DTLZ7, WFG1, WFG2, WFG3 and WFG9  for IBEA, MOEAD, NSGAIII, NSGAII and SPEA.}
		\label{tab:IGD}
		\begin{tabular}{lccccc}
			\hline\noalign{\smallskip}
			\\
			\textbf{Problem set} &&&\textbf{IGD}&&\\\cline{2-6}
			&IBEA  &MOEAD &NSGAIII&NSGAII&SPEA2\\
			\noalign{\smallskip}\hline\noalign{\smallskip}
			\textbf{DTLZ1}
			&0.874&0.878&0.874&0.876&0.891\\
			\textbf{DTLZ2}
			&0.871  &0.871 & 0.871&0.872&0.872\\
			\textbf{DTLZ3}
			&0.435 &0.404 &0.471&0.401&0.476 \\
			\textbf{DTLZ7}
			&0.884&0.876&0.880&0.880&0.877\\
			\textbf{WFG1}
			&0.871  &0.872 & 0.872&0.873&0.872\\
			\textbf{WFG2}
			&0.875  &0.875 &0.875&0.876&0.876\\
			\textbf{WFG3}
			&0.883  &0.885 &0.884&0.885&0.885\\
			\textbf{WFG9}
			&0.878  &0.878 &0.880&0.880&0.880\\
			\noalign{\smallskip}\hline
		\end{tabular}
	\end{table}
	
	The MIP DoM quality measure is a binary indicator. In Table \ref{tab:Dom_dtlz1}, for example, it is possible to see all the comparisons among the algorithms for the \textit{DTLZ} family problem set. $P$ was the solution set generated by the algorithm which was trying to dominate, and $Q$ was the solution set which was being dominated.
	
	Using MIP DoM as a unary indicator is still feasible, a straightforward idea was to merely apply a summation by the values for each algorithm in the Table \ref{tab:Dom_dtlz1}  (summing up the elements of a row, for example) .

	\begin{table}[!htb]
		\caption{HV quality indicator for problem sets generated using DTLZ1, DTLZ2, DTLZ3, DTLZ7, WFG1, WFG2, WFG3 and WFG9  for IBEA, MOEAD, NSGAIII, NSGAII and SPEA.}
		\label{tab:HV}
		\begin{tabular}{lccccc}
			\hline\noalign{\smallskip}
			\\
			\textbf{Problem set} &&&\textbf{HV}&&\\\cline{2-6}
			&IBEA  &MOEAD &NSGAIII&NSGAII&SPEA2\\
			\noalign{\smallskip}\hline\noalign{\smallskip}
			\textbf{DTLZ1}
			&1.290E+05&1.290E+05 &1.290E+05&1.290E+05&1.289E+05\\
			\textbf{DTLZ2}
			&0.709&0.667& 0.693&0.6305&0.649\\
			\textbf{DTLZ3}
			&2.469E+08  &2.469E+08 	&2.468E+08&2.469E+08 &2.467E+08\\
			\textbf{DTLZ7}
			&3.407  &3.431 &3.563&3.505&3.417\\
			\textbf{WFG1}
			&3.493  & 3.191 &3.141&2.849& 2.583\\
			\textbf{WFG2}
			&43.823  &40.427 &42.773&41.062&41.766\\
			\textbf{WFG3}
			&21.739  &19.532 &20.191&20.623&17.355\\
			\textbf{WFG9}
			&21.950  &18.254 &18.076&18.248&19.618\\
			\noalign{\smallskip}\hline
		\end{tabular}
	\end{table}
	
	For the DTLZ1 test set, which was detailed in Tables  \ref{tab:IGD} and \ref{tab:HV}, algorithms that present the best IGD were IBEA and NSGAIII. For the HV indicator, it was difficult to compare the algorithms due to inflated values. Using MIP DoM approach and the comparison among algorithms, the IBEA and NSGAIII were both indicated as the best solutions. The results are presented in Table \ref{tab:Dom_dtlz1}. For example, if the MIP DoM values were summed up, the IBEA shows a 2.424 and NSGAIII 4.175.  Moreover, when IBEA tried to dominate NSGAIII, it presented a lower value, 1.078, compared to 1.359 when NSGAIII tried to dominate IBEA.
	
	In the DTLZ2 case, IGD  indicated IBEA, MOEAD, and NSGAIII algorithms (see Table  \ref{tab:IGD}). Considering the HV, IBEA algorithm was the best one (Table  \ref{tab:HV}). Again, looking at Table \ref{tab:Dom_dtlz1}, the best values also indicated IBEA, MOEAD, and NSGAIII algorithms. However, when making two by two comparisons, MOEAD presented better value than IBEA and NSGAIII.

	\begin{table}[!htb]
		\caption{MIP DoM values for the problem set of the \textit{DTLZ} family for comparison among IBEA, MOEAD, NSGAIII, NSGAII, and SPEA2 algorithms. It must be noted that \textit{\textbf{P}} was the solution set generated by the algorithm that was trying to dominate, and \textit{\textbf{Q} } was the solution set generated by algorithm being dominated. Each solution set had \textit{NP} = \textit{NQ} = 50 solutions and \textit{M} = 3  dimensions.}
		\label{tab:Dom_dtlz1}
		\begin{tabular}{llccccc}
			\hline\noalign{\smallskip}
			&&&\textbf{DoM(P,Q)}&\\
			&&&&\textbf{Q} &\\\cline{3-7}
			\textit{\textbf{Problem set} }&\textit{\textbf{P} }&IBEA  &MOEAD &NSGAIII&NSGAII&SPEA2\\
			\noalign{\smallskip}\hline\noalign{\smallskip}
			\textbf{DTLZ1 }
			&\textit{IBEA} &0.000&0.876&1.078&0.470&0.000\\
			&\textit{MOEAD}&2.122&0.000&2.105&2.010&0.004\\
			&\textit{NSGAIII}&1.359&1.400&0.000&1.357&0.060\\
			&\textit{NSGAII}&1.759&1.664&1.759&0.000&1.759\\
			&\textit{SPEA2}&4.159&4.161&4.160&4.161&0.000\\\cline{2-7}
			\textbf{DTLZ2}
			&\textit{IBEA} &0.000&0.988&0.968&0.993&0.984\\
			&\textit{MOEAD}&0.885&0.000&0.832&0.945&0.841\\
			&\textit{NSGAIII}&0.980&0.927&0.000&0.979&0.929\\
			&\textit{NSGAII}&1.013&1.006&0.999&0.000&0.948\\
			&\textit{SPEA2}&1.060&1.060&1.060&1.071&0.000\\\cline{2-7}
			\textbf{DTLZ3}
			&\textit{IBEA} &0.000&6.920&0.008&0.008&0.006\\
			&\textit{MOEAD}&6.526&0.000&0.007&3.415&0.001\\
			&\textit{NSGAIII}&40.748&40.748&0.000&40.881&6.891\\
			&\textit{NSGAII}&14.790&20.328&22.884&0.000&0.000\\
			&\textit{SPEA2}&49.151&49.151&22.290&49.151&0.000\\\cline{2-7}
			\textbf{DTLZ7}
			&\textit{IBEA} &0.000&0.000&0.083&0.051&0.000\\
			&\textit{MOEAD}&1.716&0.000&1.745&1.716&1.701\\
			&\textit{NSGAIII}&1.012&0.065&0.000&1.012&0.032\\
			&\textit{NSGAII}&1.048&0.143&1.048&0.000&0.116\\
			&\textit{SPEA2}&1.747&1.498&1.747&1.747&0.000\\
			\noalign{\smallskip}\hline
		\end{tabular}
	\end{table}

	The result for DTLZ3 was presented in Tables  \ref{tab:IGD} and \ref{tab:HV}: for IGD, the best algorithms were NSGAII and MOEAD; and IBEA, MOEAD, and NSGAII had the best HV values. Summing up the MIP DoM values among all algorithms (as shown in Table \ref{tab:Dom_dtlz1}), the best algorithms were IBEA and MOEAD. In this case,MIP  DoM was in agreement with HV. Comparing IBEA and MOEAD directly, MOEAD presented a lower DoM.

	\begin{table}[!htb]
		\caption{MIP DoM values for the problem set of the  \textit{WFG} family for comparison among IBEA, MOEAD, NSGAIII, NSGAII, and SPEA2 algorithms. It must be noted that \textit{\textbf{P}} was the solution set generated by the algorithm that was trying to dominate, and \textit{\textbf{Q} } was the solution set generated by algorithm being dominated. Each solution set had \textit{NP} = \textit{NQ} = 50 solutions and \textit{M} = 3  dimensions.}
		\label{tab:Dom_wfg}
		\begin{tabular}{llccccc}
			\hline\noalign{\smallskip}
			&&&\textbf{DoM(P,Q)\textit{}}&\\
			&&&&\textbf{Q} &\\\cline{3-7}
			\textit{\textbf{Problem set} }&\textit{\textbf{P} }&IBEA  &MOEAD &NSGAIII&NSGAII&SPEA2\\
			\noalign{\smallskip}\hline\noalign{\smallskip}
			\textbf{WFG1 }
			&\textit{IBEA} &0.000&0.916&0.826&0.842&0.916\\
			&\textit{MOEAD}&1.046&0.000&1.007&1.007&1.007\\
			&\textit{NSGAIII}&0.989&0.827&0.000&0.941&0.989\\
			&\textit{NSGAII}&1.021&0.981&1.021&0.000&1.021\\
			&\textit{SPEA2}&0.246&0.133&0.134&0.098&0.000\\
			\textbf{WFG2}
			&\textit{IBEA} &0.000&0.774&1.092&0.860&1.078\\
			&\textit{MOEAD}&1.544&0.000&1.544&1.544&1.544\\
			&\textit{NSGAIII}&1.536&1.548&0.000&1.536& 1.497\\
			&\textit{NSGAII}&1.613&1.548&1.613&0.000&1.562\\
			&\textit{SPEA2}&1.679&1.635&1.679&1.679&0.000\\
			\textbf{WFG3}
			&\textit{IBEA} &0.000&2.997&2.997&2.997&2.998\\
			&\textit{MOEAD}&2.488&0.000&2.260&2.932&2.543\\
			&\textit{NSGAIII}&2.144&1.945&0.000&2.865&2.946\\
			&\textit{NSGAII}&1.242&1.705&1.480&0.000&3.071\\
			&\textit{SPEA2}&0.754&0.541&0.513&1.076&0.000\\
			\textbf{WFG9}
			&\textit{IBEA} &0.000&1.932&2.115&2.115&2.199\\
			&\textit{MOEAD}&2.009&0.000&2.060&2.009&2.009\\
			&\textit{NSGAIII}&2.085&1.815&0.000&2.010&2.085\\
			&\textit{NSGAII}&2.310&2.310&2.310&0.000&2.310\\
			&\textit{SPEA2}&2.219&2.219&2.219&2.219&0.000\\
			\noalign{\smallskip}\hline
		\end{tabular}
	\end{table}

	In DTLZ7 problem set, the best HV values were NSGAIII and NSGAII. Considering IGD, the best values were for MOEAD and SPEA2. Applying the same summing up approach as before, and using Table \ref{tab:Dom_dtlz1}, the result was that NSGAIII and NSGAII generated the best candidate solutions. Again, comparing NSGAIII and NSGAII, the best one was NSGAIII with a DoM(P, Q) = 1.012.
	
	The same experiment was done for the \textit{WFG} family. For WFG1 in Tables \ref{tab:IGD} and \ref{tab:HV}, the best algorithm was IBEA for IGD, and HV indicator. In Table \ref{tab:Dom_wfg}, the best algorithm was SPEA2, and the second one was IBEA. However, comparing SPEA2 with IBEA, SPEA2 had a lower value of MIP DoM, DoM(SPEA2, IBEA) = 0.246, in contrast with  DoM(IBEA, SPEA2) = 0.916. This behavior could be observed in Figure \ref{fig:WFG9_sample}. Visually, it is possible to see a similarity between the two solution sets. However, it is relevant to note the graphics' scale in the axis. The SPEA2 solution set appeared more convergent than IBEA (the spread are similar, but SPEA2 presented minor values in the axis).  
	
	In WFG2, the best values for IGD were entirely tied, with a little difference (see Table \ref{tab:IGD}). Considering HV, the best one was the IBEA algorithm (Table \ref{tab:HV}). Using DoM, detailed in Table \ref{tab:Dom_wfg}, it was possible to see the same characteristic as reported by IGD. Nonetheless, if the algorithms were compared two by two, it was possible to indicate IBEA as the algorithm which presented the lower MIP DoM values.   
	
	Using WFG3, it was possible to observe little differences for the IGD indicator among the other algorithms. For HV, there was an indication that IBEA had a greater value, but with other values next to it. Assessing MIP DoM in Table \ref{tab:Dom_wfg}, there was an indication that SPEA2 had better values. Furthermore, it was still possible to analyze the sets graphically in Figure \ref{fig:WFG9_sample}. A visual and qualitative analysis comparing IBEA with SPEA2 graphics could show how DoM acknowledged the spread's solution set. All the solution sets had a similar spread and cardinality (it is possible to see a kind of straight line formed by points), and SPEA2 was one notable exception.  
	
	Finally, for the WFG9 problem set, Tables \ref{tab:IGD} and \ref{tab:HV} showed that for IGD,  algorithms IBEA and MOEAD had lower values, but the values were next to each other. For HV, the IBEA algorithm had a better value. Considering DoM, presented in Table \ref{tab:Dom_wfg}, all three algorithms were competitive: IBEA, MOEAD, and NSGAIII. Comparing IBEA and MOEAD, for example, the best MIP DoM value was found for IBEA: DoM(IBEA,MOEAD) = 1.932.
	
	The discrepancy for the WFG1 and WFG3 experiments, among DoM, IGD, and HV was something that drew attention. For DoM, the optimal value was obtained using the model described in \ref{mip_dom_model}. However, the reference sets for IGD and HV were generated in a simple manner. This might be improved in assessing the convergence among the indicators. Another possibility would be to increase the number of solutions in each set.

	\begin{figure}[!htb]
		\centering
		{\includegraphics[width=\linewidth]{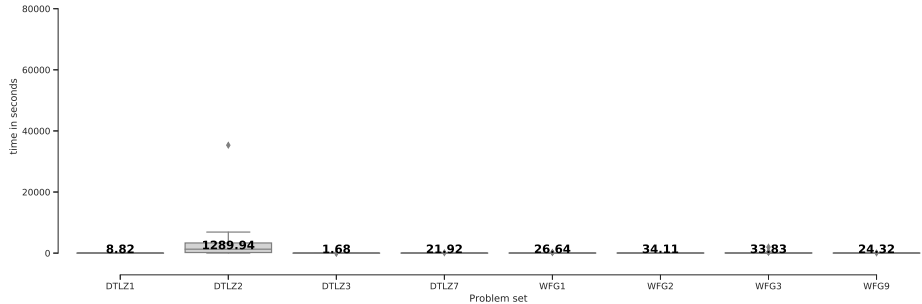} }
		\caption[Time spent in seconds]{A box plot with the time spent by the GUROBI solver \cite{gurobi} to solve the MIP DoM model for the DTLZ1, DTLZ2, DTLZ3, DTLZ7, WFG1, WFG2, WFG3, and WFG9 problem sets with \textit{NP} = \textit{NQ} = 50 points and \textit{M} = 3 objectives. }
		\label{fig:BP_time_50_3}%
	\end{figure}

	The GUROBI solver \cite{gurobi}  used branch and bound,  branch and cut, and other solver capabilities, such as: Gomory cuts, flow cover, clique, and others. The minimum gap established was five percent, and the only additional parameter used was MIPFocus, which helps to improve the best bound during the execution. For some models, the time spent was small. Figure \ref{fig:BP_time_50_3} shows a box plot for the results with the medians highlighted in bold. All the models were solved, and the median was less than 35 seconds. The only exception was DTLZ2, which presented the highest median time spent, $\sim1289$. Furthermore, there were two cases in high discrepancy, more specifically: IBEA trying to dominate NSGAII, with $\sim89444$ seconds, and NSGAIII trying to dominate NSGAII with $\sim35325$ seconds. These two cases had some common characteristics of what the model had spent more time than others (high density between the solution sets, for example).
	
\begin{figure}[!htb]
		\centering
		{\includegraphics[width=\linewidth]{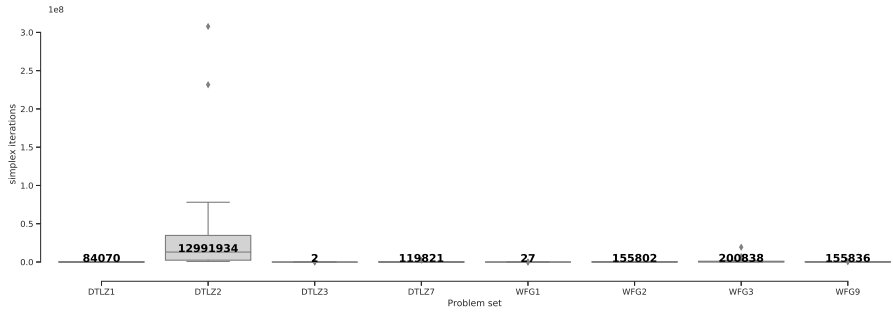} }
		\caption[Simplex iterations]{A box plot with the simplex iterations used by the GUROBI solver \cite{gurobi} to solve the DoM model for the DTLZ1, DTLZ2, DTLZ3, DTLZ7, WFG1, WFG2, WFG3, and WFG9 problem sets with \textit{NP} = \textit{NQ} = 50 points and \textit{M} = 3 objectives.  }
		\label{fig:BP_simplex_iterations_50_3}%
	\end{figure}

	In Figure \ref{fig:BP_simplex_iterations_50_3}, a piece of similar information was depicted: the number of simplex iterations by the solver. Again, in some cases, such as DTLZ3 and WFG1, the number of iterations was small. The same discrepancy cases, as shown in Figure \ref{fig:BP_time_50_3}, happened here with DTLZ2. 
	
	\subsubsection{Many-Objective cases}
	
	In this Section, the goal was to verify if the MIP DoM approach could be applied in many-objective scenario. The same problem sets and algorithms from the last experiment were used. It worthy to note that the number of points was established based on some works (\cite{8027123} and \cite{Yang2019}, for example).
	
	In Table \ref{tab:Dom_dtlz_100_5}, problem sets from the \textit{DTLZ} family are presented. At this time, the problems had $M = 5$ objectives and $NP = NQ = 100$ solutions. The non-dominated points were generated using IBEA, MOEAD, NSGAIII, NSGAII, and SPEA2. For the DTLZ1 problem set, the best algorithms were IBEA and SPEA2, but when they were compared, DoM(IBEA, SPEA2) = 0.525 and DoM(SPEA2, IBEA) = 2.278, IBEA was still better than SPEA2. 
	
	Following this idea, for DTLZ2, the IBEA and SPEA2 were more competitive than the others. In this case, two by two comparison using DoM was very tight (1.000 and 1.060, for example). However, IBEA presented a lower value. Considering DTLZ3, again, IBEA and SPEA2 were the best ones, and IBEA exhibited the best value, compared with SPEA2 or any other algorithm. Finally, for DTLZ7, it was possible to observe a curious fact. In this case, SPEA2 already weakly dominated all other solution sets.
	
	\begin{table}[!htb]
		\caption{MIP DoM values for the problem set of the \textit{DTLZ} family for comparison among IBEA, MOEAD, NSGAIII, NSGAII, and SPEA2 algorithms. It must be noted that \textit{\textbf{P}} was the solution set generated by the algorithm that was trying to dominate, and \textit{\textbf{Q} } was the solution set generated by algorithm being dominated. Each solution set had \textit{NP} = \textit{NQ} = 100 solutions and \textit{M} = 5  dimensions.}
		\label{tab:Dom_dtlz_100_5}
		\begin{tabular}{llccccc}
			\hline\noalign{\smallskip}
			&&&\textbf{DoM(P,Q)\textit{}}&\\
			&&&&\textbf{Q} &\\\cline{3-7}
			\textit{\textbf{Problem set} }&\textit{\textbf{P} }&IBEA  &MOEAD &NSGAIII&NSGAII&SPEA2\\
			\noalign{\smallskip}\hline\noalign{\smallskip}
			\textbf{DTLZ1 }
			&\textit{IBEA} &0.000&0.000&0.000&0.522&0.525\\
			&\textit{MOEAD}&40.616&0.000&4.061&40.616&40.616\\
			&\textit{NSGAIII}&14.749&14.749&0.000&14.749&14.749\\
			&\textit{NSGAII}&17.595&17.595&17.595&0.000&17.595\\
			&\textit{SPEA2}&2.278&0.000&0.000&0.882&0.000\\\cline{2-7}
			\textbf{DTLZ2}
			&\textit{IBEA} &0.000&1.001&0.357&0.476&1.000\\
			&\textit{MOEAD}&1.250&0.000&1.294&1.269&1.299\\
			&\textit{NSGAIII}&1.147&1.184&0.000&1.147&1.203\\
			&\textit{NSGAII}&1.074&1.074&1.074&0.000&1.074\\
			&\textit{SPEA2}&1.060&0.457&0.749&1.060&0.000\\\cline{2-7}
			\textbf{DTLZ3}
			&\textit{IBEA} &0.000&16.231&16.231&0.276&0.333\\
			&\textit{MOEAD}&450.66&0.000&450.66&450.66&450.66\\
			&\textit{NSGAIII}&488.037&488.037&0.000&488.037&488.037\\
			&\textit{NSGAII}&356.921&170.168&350.819&0.000&356.95\\
			&\textit{SPEA2}&59.66&0.000&0.000&2.055&0.000\\\cline{2-7}
			\textbf{DTLZ7}
			&\textit{IBEA} &0.000&0.733&1.677&1.597&1.710\\
			&\textit{MOEAD}&2.848&0.000&2.811&2.834&2.811\\
			&\textit{NSGAIII}&0.729&0.065&0.000&0.543&0.032\\
			&\textit{NSGAII}&2.707&2.336&2.681&0.000&2.707\\
			&\textit{SPEA2}&0.000&0.000&0.000&0.000&0.000\\
			\noalign{\smallskip}\hline
		\end{tabular}
	\end{table}
	
	In Table \ref{tab:Dom_WFG_100_5},  problems from the WFG family were tested, with the same proposal as before. For WFG1, the DoM value comparison was tight for IBEA, MOEAD, NSGAIII, and SPEA2. The ranking, using the idea of move dominance summation, was for NSGAIII and MOEAD as the best ones.
	
	For the WFG2 problem, all the lower values reference IBEA algorithm as the best solution set for this problem. Considering WFG3, there was, again, a curious fact: the solution set generated by IBEA almost dominated all other results, the only exception was for MOEAD. Finally, for WFG9, it was clear that the solution set generated by IBEA was lower than the others.
	
	Another relevant fact which is important to observe: some MIP DoM values were repeated, when a solution $P$ was trying to dominate $Q$, for example. See in Table \ref{tab:Dom_WFG_100_5}, the problem set WFG1 and $P$ as NSGAII, or the problem set WFG2 and $P$ as MOEAD, or even the problem set WFG9 and $P$ as SPEA2. These cases raised one question: What are the points in $P$ that were updated to dominate $Q$? In all the last cases, the same points $P'$ were generated by the same $P$. For example, for the case with the solution set WFG1 and $P$ as NSGAII, $p_{8}$ was always the point where $p'_{8}$ was generated, and this was the final solution for all the cases (i.e., in each case the values updated in each objective were different, but all the final solutions used $p_{8}$). The same fact happened with for all last cited cases.
	
\begin{figure}[!htb]
		\centering
		{\includegraphics[width=12cm]{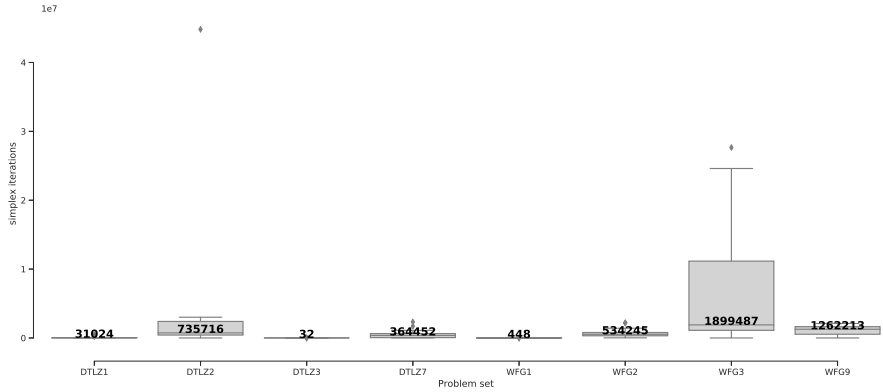} }
		\caption[Simplex iterations]{A box plot with the simplex iterations used by the GUROBI solver \cite{gurobi} to solve the DoM model for the DTLZ1, DTLZ2, DTLZ3, DTLZ7, WFG1, WFG2, WFG3, and WFG9 problem sets with \textit{NP} = \textit{NQ} = 100 points and \textit{M} = 5 objectives.  }
		\label{fig:BP_simplex_iterations_100_5}%
	\end{figure}

	In Figure \ref{fig:BP_simplex_iterations_100_5}, the number of simplex iterations was presented. There was a discrepancy in DTLZ2, which demanded much more iterations than the others. This was when IBEA was trying to dominate SPEA2. The time spent in this extreme case was $\sim10245$  seconds.

	\begin{table}[!htb]
		\caption{MIP DoM values for the problem set of the \textit{WFG} family for comparison among IBEA, MOEAD, NSGAIII, NSGAII, and SPEA2 algorithms. It must be noted that \textit{\textbf{P}} was the solution set generated by the algorithm that was trying to dominate, and \textit{\textbf{Q} } was the solution set generated by algorithm being dominated. Each solution set had \textit{NP} = \textit{NQ} = 100 solutions and \textit{M} = 5  dimensions.}
		\label{tab:Dom_WFG_100_5}
		\begin{tabular}{llccccc}
			\hline\noalign{\smallskip}
			&&&\textbf{DoM(P,Q)\textit{}}&\\
			&&&&\textbf{Q} &\\\cline{3-7}
			\textit{\textbf{Problem set} }&\textit{\textbf{P} }&IBEA  &MOEAD &NSGAIII&NSGAII&SPEA2\\
			\noalign{\smallskip}\hline\noalign{\smallskip}
			\textbf{WFG1 }
			&\textit{IBEA} &0.000&0.398&0.398&0.235&0.396\\
			&\textit{MOEAD} &0.419&0.000& 0.401&0.062&0.422\\
			&\textit{NSGAIII} &0.391&0.390&0.000&0.129&0.391\\
			&\textit{NSGAII} &1.133&1.133&1.133&0.000&1.133\\
			&\textit{SPEA2} &0.426&0.426&0.426&0.229&0.000\\\cline{2-7}
			\textbf{WFG2}
			&\textit{IBEA} &0.000&1.272&1.601&1.601&1.655\\
			&\textit{MOEAD} &2.457&0.000&2.457&2.457&2.457\\
			&\textit{NSGAIII} &2.345&2.345&0.000&2.345&2.345\\
			&\textit{NSGAII} &1.873&1.873&1.873&0.000&1.873\\
			&\textit{SPEA2} &1.708&1.697&1.773&1.708&0.000\\\cline{2-7}
			\textbf{WFG3}
			&\textit{IBEA} &0.000&0.038&0.000&0.000&0.000\\
			&\textit{MOEAD} &5.390&0.000&4.361&4.191&4.692\\
			&\textit{NSGAIII} &6.460&5.914&0.000&4.997&5.651\\
			&\textit{NSGAII} &6.262&6.262&6.262&0.000&6.262\\
			&\textit{SPEA2} &5.901&5.368&5.207&4.571&0.000\\\cline{2-7}
			\textbf{WFG9}
			&\textit{IBEA} &0.000&2.152&2.151&2.151&2.166\\
			&\textit{MOEAD} &4.346&0.000&4.300&4.237&4.346\\
			&\textit{NSGAIII} &5.515&5.515&0.000&5.515&5.515\\
			&\textit{NSGAII} &3.188&3.189&3.189&0.000&3.189\\
			&\textit{SPEA2} &4.427&4.427&4.427&4.427&0.000\\
			\noalign{\smallskip}\hline
		\end{tabular}
	\end{table}

	The final experiment was made to verify if the MIP DoM approach was viable in a higher dimensional objective space . The test was done considering \textit{M} = 10 objectives and \textit{NP} = \textit{NQ} = 200. It was necessary to reduce the number of algorithms. In Platypus \cite{Brockhoff:2019:Platypus}, there was a parameter to specify the population size; however, for MOEAD and NSGAIII, this parameter was ignored for ten dimensions. For this reason, in Table \ref{tab:Dom_dtlz_200_10}, these algorithms were not presented. 
	
	Table \ref{tab:Dom_dtlz_200_10} shows that the MIP DoM approach could be computed using these test sets. There were some cases in which the model spent 2.05 seconds, with 2556 simplex iterations for one case, or 29.97 seconds, with 45674 simplex iterations, for other cases. These two examples were the fastest ones. On the other hand, there was another one, in which 180468.60 seconds was spent, with 107009076 simplex iterations, presenting the worst case. 
	
	For the fastest cases, there were many points in the solution set $Q$, that already were dominated by some point in $P$. For these cases, the model and the solver exploited this characteristic, and they were capable of finding the best solution in a reasonable time.
	
	On the other hand, the worst case spent almost 51 hours. In this situation, there was no a priori dominance among the points in each solution set.  One thing that can be better investigated is the density between the solution sets when there is a chance for $p$ to change its position to dominate some $q's$, considering a $g$  group. In areas with high density, it is plausible that there are many options concerning $P$ and $Q$, and the branch and bound dynamic become more challenging.  
	
	One last interesting fact observed in the solutions generated by the MIP DoM model is related to the number of $p's$ that was altered to obtain the MIP DoM value. For this final  experiment, for example, in which the space was large, the maximum number of $p's$ that had their objectives updated was six. The minimum number of $p's$ was one. For these last cases, the MIP DoM recommendation was to change only one point, and that point was capable of dominating all the solutions and generating the best MIP DoM value.

	\begin{table}[!htb]
		\caption{MIP DoM values for \textit{DTLZ} and \textit{WFG} family problem set for comparisons among IBEA, NSGAII, and SPEA2 algorithms.  It must be noted that \textit{\textbf{P}} was the solution set generated by the algorithm that was trying to dominate, and \textit{\textbf{Q} } was the solution set generated by algorithm being dominated. Each solution set had \textit{NP} = \textit{NQ} = 200 solutions and \textit{M} = 10  dimensions.}
		\label{tab:Dom_dtlz_200_10}
		\begin{tabular}{llcccc}
			\hline\noalign{\smallskip}
			&\textbf{DoM(P,Q)\textit{}}&\\
			&&&\textbf{Q} \\\cline{3-5}
			\textit{\textbf{Problem set} }&\textit{\textbf{P} }&IBEA  &NSGAII&SPEA2\\
			\noalign{\smallskip}\hline\noalign{\smallskip}
			\textbf{DTLZ1 }
			&\textit{IBEA} &0.000&0.926&1.102\\
			&\textit{NSGAII} &52.666&0.000&52.666\\
			&\textit{SPEA2} &9.585&4.986&0.000\\\cline{2-5}
			\textbf{DTLZ2}
			&\textit{IBEA} &0.000&0.981&1.024\\
			&\textit{NSGAII} &1.420&0.000&1.430\\
			&\textit{SPEA2} &1.843&1.843&0.000\\\cline{2-5}
			\textbf{DTLZ3}
			&\textit{IBEA} &0.000&0.004&13.465\\
			&\textit{NSGAII} &627.13&0.000&598.422\\
			&\textit{SPEA2} &425.427&425.427&0.000\\\cline{2-5}
			\textbf{DTLZ7}
			&\textit{IBEA} &0.000&0.365&0.181\\
			&\textit{NSGAII} &3.305&0.000&3.305\\
			&\textit{SPEA2} &3.757&3.757&0.000\\\cline{2-5}
			\textbf{WFG1}
			&\textit{IBEA} &0.000&0.288&0.416\\
			&\textit{NSGAII} &1.261&0.000&1.261\\
			&\textit{SPEA2} &0.832&0.832&0.000\\\cline{2-5}
			\textbf{WFG2}
			&\textit{IBEA} &0.000&5.653&5.777\\
			&\textit{NSGAII} &3.568&0.000&3.523\\
			&\textit{SPEA2} &3.530&3.530&0.000\\\cline{2-5}
			\textbf{WFG3}
			&\textit{IBEA} &0.000&0.000&0.000\\
			&\textit{NSGAII} &14.487&0.000&14.487\\
			&\textit{SPEA2} &14.231&13.358&0.000\\\cline{2-5}
			\textbf{WFG9}
			&\textit{IBEA} &0.000&15.706&15.706\\
			&\textit{NSGAII} &19.370&0.000&19.404\\
			&\textit{SPEA2} &15.949&15.949&0.000\\
			\noalign{\smallskip}\hline
		\end{tabular}
	\end{table}

	As can be viewed in Equation (\ref{comb_plinha}), the combinatorial nature of the MIP DoM calculation is one reason which has made it difficult to use. The computational time increases exponentially with the size of the solution set. The mixed integer programming is NP hard in general \cite{DBLP:books/daglib/0023873}. Despite this fact, the MIP approach has managed to be viable in three objectives with 50 solutions, five objectives with 100 solutions, and ten objectives with 200 solutions for some DTLZ and WFG family test cases. 
	
	Concerning the MIP model, some issues need to be mentioned:
	
	\begin{enumerate}[i)]
		\item In the LP relaxation for  $xp_{(i)}$ and $xpq{(i,j)}$ (in the first two constraints in Equation (\ref{restriction_xpodB})), the integrality can not be ensured; this problem is related in \cite{articleKlotz} causing, in some cases, a slow advance in the best bound;
		\item The last problem could be solved if we knew, a priori, the number of $p's$ that would dominate $Q$. For the last experiment, for example, it was observed that in a high number of dimensions, a low percentage of points was used to dominate the other set. We believe that imposing a constraint that limits the number of $p's$ can improve the number of simplex iterations.
		\item The coefficient matrix to bound $p^{'}_i$ presented a large amplitude in some cases, and it causes numerical difficult to solvers. 
		
	\end{enumerate}

	\section{Conclusion}\label{conclusion}
	
	The DoM binary indicator considers the minimum move of one set to weakly dominate the other set. Given two solution sets $P$ and $Q$, the DoM of $P$ to $Q$ is the minimum total distance of moving some points in $P$ in way that any point in $Q$ is weakly dominated by at least one point in $P$. The indicator is Pareto compliant and does not demand any parameters or reference sets. The handicap about DoM is its combinatorial calculation nature making it difficult to be applied to problems with three or more objective functions.
	
	We explored a MIP approach calculation of DoM, trying to understand if it was a viable way to deal with the problem. The idea was to use the polyhedral method, as it has been proved to be successful in practice with other issues. A comparison with artificial bidimensional examples, the same proposed by the  DoM's authors was made, with the same results. Additionally, some problems in three dimensions showed that MIP DoM results were in agreement and compliant compared with other common indicators used in literature (IGD and HV). Finally, experiments using five and ten dimensions were executed and had shown that MIP DoM model was viable and, in some cases with reasonable computational time. All the experiments have used a fixed number of solutions in each set (i.e., $NP = NQ = \{50, 100, 200\}$, for $M = \{3, 5, 10\}$ dimensions). To the best of our knowledge, even with these limitations, an exact method to calculate DoM in three or more dimensions was not known until now.
	
	For some experiments, high variability was found in the time spent by the model resolution. It was observed that this fact was inherent to the distribution or density, involving the two solution sets, $P$ and $Q$. This issue is something that deserves to be best investigated in order to improve the model, creating, for example, new constraints in the model.
	
	In the MIP approach for DoM, other directions also deserve to be investigated: i) How to efficiently calculate DoM using MIP and exploiting some inherent solution set features; ii) How to define a priori the minimum number of $p^{'}_i$ to be moved in an attempt to dominate $Q$ and what benefits it can bring to the MIP model. What the relation between the distributions of $P$ and $Q$, and the minimum number of $p^{'}_i$ to be moved is; iii) How to improve the amplitude of the coefficients' matrix  in the MIP model, and if this is this capable of improving the convergence model time. Such questions were posed here to emphasize possible future research directions. The tests with more dimensions and experiments with other problems test sets are something that must be done, as well.
	
	Finally, MIP DoM is an indicator that can represent all the solution sets quality facets. Moreover, it proposes a more natural and intuitive relation when comparing solution sets using a measure compatible with the Pareto dominance concept, mainly in high dimensional objective space. The MIP DoM, like its predecessor, does not require prior problem knowledge and any other additional parameter. With MIP DoM Model, a viable method to calculate the DoM indicator was proposed to try to incentive its general use for compare solution sets or even to guide some algorithms in the search for solutions with useful quality features.
	

	%
	%

	\bibliographystyle{spmpsci}      
	\bibliography{jogodom}   
	

\end{document}